\documentclass[12pt]{amsart}
\usepackage{amscd,verbatim}
\usepackage{amssymb}
\usepackage[all]{xy}
\usepackage[colorlinks,linkcolor=blue,citecolor=blue,urlcolor=red]{hyperref}

\newcommand{\un}{\mathbf{1}}
\newcommand{\A}{\mathbf{A}}

\newcommand{\G}{\mathbf{G}}

\renewcommand{\P}{\mathbf{P}}

\newcommand{\Z}{\mathbb{Z}}

\newcommand{\sC}{\mathcal{C}}

\newcommand{\sO}{\mathcal{O}}

\newcommand{\bZ}{\mathbb{Z}}

\newcommand{\Xb}{{\overline{X}}}

\newcommand{\Cor}{\operatorname{\mathbf{Cor}}}

\newcommand{\HI}{{\operatorname{\mathbf{HI}}}}

\newcommand{\Rec}{{\operatorname{\mathbf{Rec}}}}

\newcommand{\PST}{{\operatorname{\mathbf{PST}}}}

\newcommand{\NST}{\operatorname{\mathbf{NST}}}

\newcommand{\DM}{\operatorname{\mathbf{DM}}}

\newcommand{\MDM}{\operatorname{\mathbf{MDM}}}

\newcommand{\Hom}{\operatorname{Hom}}
\newcommand{\uHom}{\operatorname{\underline{Hom}}}

\newcommand{\Ker}{\operatorname{Ker}}

\newcommand{\Coker}{\operatorname{Coker}}

\newcommand{\Pic}{\operatorname{Pic}}

\newcommand{\Spec}{\operatorname{Spec}}

\newcommand{\Sm}{\operatorname{\mathbf{Sm}}}
\newcommand{\Sch}{\operatorname{\mathbf{Sch}}}
\newcommand{\Ab}{\operatorname{\mathbf{Ab}}}

\newcommand{\by}{\xrightarrow}

\newcommand{\iso}{\by{\sim}}

\newcommand{\rec}{{\operatorname{rec}}}

\newcommand{\tr}{{\operatorname{tr}}}

\newcommand{\eff}{{\operatorname{eff}}}

\renewcommand{\o}{{\operatorname{o}}}
\newcommand{\op}{{\operatorname{op}}}

\newcommand{\Zar}{{\operatorname{Zar}}}
\newcommand{\Nis}{{\operatorname{Nis}}}

\newcommand{\inj}{\hookrightarrow}
\newcommand{\Inj}{\lhook\joinrel\longrightarrow}

\newcommand{\id}{{\operatorname{Id}}}

\newcommand{\pd}{{\partial}}

\newcommand{\ch}{{\operatorname{ch}}}

\newcommand{\CH}{{\operatorname{CH}}}

\renewcommand{\lim}{\operatornamewithlimits{\varprojlim}}
\newcommand{\colim}{\operatornamewithlimits{\varinjlim}}

\newcommand{\ol}{\overline}

\newcommand{\car}{\operatorname{char}}

\renewcommand{\phi}{\varphi}
\renewcommand{\epsilon}{\varepsilon}
\renewcommand{\div}{\operatorname{div}}

\newcommand{\MNST}{\operatorname{\mathbf{MNST}}}

\newcommand{\MCor}{\operatorname{\mathbf{MCor}}}
\newcommand{\MP}{\operatorname{\mathbf{MSm}}}

\newcommand{\MPST}{{\operatorname{\mathbf{MPST}}}}
\newcommand{\CI}{{\operatorname{\mathbf{CI}}}}
\newcommand{\SCRec}{{\operatorname{\mathbf{RSC}}}}
\newcommand{\RSC}{{\operatorname{\mathbf{RSC}}}}
\newcommand{\Bl}{{\mathbf{Bl}}}

\newcommand{\bcube}{{\ol{\square}}}
\newcommand{\cube}{\square}

\newcommand{\Mb}{{\overline{M}}}

\def\rmapo#1{\overset{#1}{\longrightarrow}}

\newcommand{\ulMCor}{\operatorname{\mathbf{\underline{M}Cor}}}

\def\bZ{\mathbb{Z}}
\def\Ztr{\bZ_\tr}

\def\Image{\mathrm{Image}}

\def\Xinf{X^\infty}

\newcommand{\Vb}{\overline{V}}

\def\Vb{\overline{V}}
\def\Wb{\overline{W}}
\def\WbNo{\Wb^{N,o}}
\def\bP{\mathbf{P}}

\def\CXY#1#2{C_{#1}(\Xb|Y)(#2)}

\def\Minf{M_\infty}
\def\qaq{\;\text{ and }\;}
\def\qwith{\;\text{ with}\;}
\def\hM#1{h_0^{\bcube}(#1)}
\def\hMM#1{h^0_{\bcube}(#1)}

\def\hMw#1{h_0(#1)}
\def\ep{\epsilon}
\def\ta{\tilde{a}}
\def\aVNis{a^V_\Nis}

\newcounter{spec}
\newenvironment{thlist}{\begin{list}{\rm{(\roman{spec})}}
{\usecounter{spec}\labelwidth=20pt\itemindent=0pt\labelsep=10pt}}
{\end{list}}

\newtheorem{lemma}{Lemma}[subsection]

\newtheorem{Th}{Theorem}

\newtheorem{thm}[lemma]{Theorem}

\newtheorem{prop}[lemma]{Proposition}
\newtheorem{proposition}[lemma]{Proposition}
\newtheorem{cor}[lemma]{Corollary}
\newtheorem{corollary}[lemma]{Corollary}

\theoremstyle{definition}

\newtheorem{definition}[lemma]{Definition}

\theoremstyle{remark}
\newtheorem{Rm}[Th]{Remark}
\newtheorem{qn}[lemma]{Question}
\newtheorem{rk}[lemma]{Remark}
\newtheorem{remark}[lemma]{Remark}
\newtheorem{remarks}[lemma]{Remarks}

\newtheorem{claim}[lemma]{Claim}

\numberwithin{equation}{section}

\begin{document}
\title{Reciprocity sheaves, II}
\author{Bruno Kahn}
\address{IMJ-PRG\\Case 247\\
4 place Jussieu\\
75252 Paris Cedex 05\\
France}
\email{bruno.kahn@imj-prg.fr}
\author{Shuji Saito}
\address{Graduate School of Mathematical Sciences\\
University of Tokyo\\
3-8-1 Komaba, \\ Tokyo 153-8941\\ Japan}
\email{sshuji@msb.biglobe.ne.jp}
\author[Takao Yamazaki]{Takao Yamazaki}
\address{Institute of Mathematics\\ Tohoku University\\ Aoba\\ Sendai 980-8578\\ Japan}
\email{takao.yamazaki.b6@tohoku.ac.jp}
\date{February 18, 2021}
\thanks{The first author acknowledges the support of Agence Nationale de la Recherche (ANR) under reference ANR-12-BL01-0005. 
The second author is supported by JSPS KAKENHI Grant (15H03606).
The third author is supported by JSPS KAKENHI Grant (15K04773). 
}

\begin{abstract}
We exhibit an intimate relationship between  ``reciprocity sheaves'' from \cite{rec} and
``modulus sheaves with transfers'' from \cite{MwM1, MwM2}.
\end{abstract}


\maketitle

\setcounter{tocdepth}{1}
\tableofcontents

\section*{Introduction}
\setcounter{subsection}{1}

This paper is a synthesis of \cite{rec} and \cite{MwM1, MwM2}; part of it uses results of \cite{rs2} and \cite{shuji}.

In \cite{rec},  we introduced reciprocity (pre)sheaves
as a generalization of 
Voevodsky's homotopy invariant (pre)sheaves with transfers, which are the main building block 
for constructing his triangulated categories of motives in \cite{voetri}. 
(From now on, we shall replace \emph{homotopy invariant} by  \emph{$\A^1$-invariant} for clarity.) Let $\Sm$ be the category of separated smooth schemes of finite type over $k$. 
There is an additive category $\Cor$ which has the same objects as $\Sm$ and whose morphisms are finite correspondences; the category $\PST$ of \emph{presheaves with transfers} is defined as the additive dual of $\Cor$  \cite[Lect. 1 and 2]{mvw}. A presheaf with transfers $F$ is \emph{$\A^1$-invariant} if the projection $X\times \A^1\to X$ induces an isomorphism $F(X) \iso F(X\times \A^1)$ for all $X\in \Sm$.
Let $\HI\subset \PST$ be the full subcategory of $\A^1$-invariant presheaves with transfers. 
The reciprocity presheaves defined in \cite{rec} form a full subcategory $\Rec\subset \PST$, which contains $\HI$. 

In this paper, we introduce a new full subcategory $\RSC\subset \PST$ which is fairly
close to $\Rec$ and fits better with the new framework of \emph{modulus presheaves
of transfers}. The latter were introduced in \cite{ksy} 
to construct a new triangulated category $\MDM^\eff$ of motivic nature which enlarges Voevodsky's triangulated category of motives $\DM^\eff$ \cite[Lect. 14]{mvw}. Due to problems encountered in \cite{ksy}, this theory was refounded in \cite{MwM1,MwM2}  and \cite{MwM3}. In this paper, we only use results from \cite{MwM1} and \cite{MwM2}, except for the tensor structure on $\MCor$.

To give an idea of how one defines $\Rec$ and $\RSC$, we need to reformulate the definition of $\A^1$-invariance. 
Recall \cite[Lem. 2.16]{mvw} that the inclusion $\HI\to \PST$ has a left adjoint \begin{equation}\label{h0A}
h_0^{\A^1}: \PST\to \HI.
\end{equation} 
Thus $F\in \PST$ is in $\HI$ if and only if for any $X\in \Sm$ and $a\in F(X)$, the map $\Ztr(X) \to F$ in $\PST$ associated to $a$ by Yoneda's lemma factors through $h_0^{\A^1}(X):=h_0^{\A^1}(\Z_\tr(X))$, where $\Ztr(X)$ is the presheaf with transfers represented by $X$.
To define reciprocity presheaves, we introduced in \cite{rec} bigger quotients $h(M)$ of $\Ztr(X)$ associated to 
a \emph{modulus pair} $M=(\Xb,\Xinf)$, consisting of a proper scheme $\Xb$ over $k$
and an effective Cartier divisor $\Xinf$ on it,
such that $X=\Xb \setminus |\Xinf|$.
Then a presheaf with transfers $F\in \PST$ belongs to $\Rec$  \cite[Definition 2.1.3]{rec} if 

\begin{quote} (*) For any quasi-affine $X\in \Sm$ and any $a\in F(X)$, the associated map $a:\Ztr(X) \to F$ factors through $h(M)$ for some $M$ as above. 
\end{quote}

The definition of the quotients $h(M)$ is very technical; it is inspired by the theorem of Rosenlicht-Serre on reciprocity for morphisms from curves to commutative algebraic groups \cite[Ch. III]{gacl}. 

Let us now recall the story of \cite{MwM1, MwM2}.
We define a category $\MCor$: 
its objects are modulus pairs $M=(\Xb,\Xinf)$ as above such that $M^\circ=\Xb-|\Xinf| \in \Sm$: this is called the \emph{interior} of $M$.
Morphisms of $\MCor$ are finite correspondences between interiors satisfying an admissibility condition with respect to $\Xinf$ (see Definition \ref{def:mcor}). 
Let $\MPST$ be the additive dual of $\MCor$. There is a pair of adjoint functors 
\begin{equation*}\label{eq;adjunction}
\MPST\begin{smallmatrix}\omega_!\\ \longrightarrow\\ \omega^*\\ \longleftarrow\\
\end{smallmatrix}\PST. 
\end{equation*}
Here $\omega^*$ is induced by the ``interior'' functor 
\[\omega:\MCor \to \Cor:\; (\Xb,\Xinf) \allowbreak\mapsto \Xb \setminus |\Xinf|,\] 
and $\omega_!$ is the left Kan extension of $\omega$.

Let $\bcube=(\P^1,\infty)\in \MCor$: we say that $F\in \MPST$ is \emph{$\bcube$-invariant}
if the ``projection'' $M\otimes \bcube\to M$ induces an isomorphism $F(M) \iso F(M\otimes \bcube)$ for all $M\in \MCor$ (see \S \ref{sect:mpst} for the monoidal structure $\otimes$ on $\MCor$).
We let $\CI\subset \MPST$ denote the full subcategory of $\bcube$-invariant objects.

We show in Theorem \ref{rem:monoidal-ci} that the inclusion $\CI\to \MPST$ has a left adjoint $h_0^{\bcube}: \MPST\to \CI$. Define $h_0(M)\in \PST$ to be $\omega_!h_0^{\bcube}\Ztr(M)$, where $\Ztr(M)\in \MPST$ is the presheaf represented by $M$ and $\omega_!$ is as before. Then $\RSC$ is the full subcategory of $\PST$ consisting of those presheaves verifying Condition (*) above, modified by dropping the quasi-affine condition on $X$ and replacing $h(M)$ by $h_0(M)$.

Our main results are the following.

\begin{Th}\label{T1} \
\begin{enumerate}
\item (Corollary \ref{c4.2}). We have $\HI\subset \RSC$.
\item (Th. \ref{cor:cube-inv-rec2} and Prop. \ref{prop:omega-ci}). 
We have $\omega_!(\CI)= \RSC$.
The induced functor $\omega_\CI: \CI\to \RSC$ 
has a fully faithful right adjoint $\omega^\CI:\RSC\to \CI$.
\end{enumerate}
\end{Th}

\begin{Th}[]\label{C1}
\label{prop:sc-rec-rec-intro}\
\begin{enumerate}
\item (Th. \ref{thm:h-h0}). Let $M=(\ol{X}, Y) \in \MCor$ be such that
$X:=\ol{X} \setminus |Y|$ is quasi-affine. Then $h_0(M)=h(M)$.
Consequently, we have $\SCRec \subset \Rec$.
\item (Cor. \ref{prop:sc-rec-rec}).
We have 
\[\RSC\cap \NST = \Rec\cap\NST.\]
Here, $\NST\subset \PST$ is the full subcategory of 
\emph{Nisnevich sheaves with transfers} \cite{voetri}. 
\end{enumerate}
\end{Th}

\enlargethispage*{20pt}

Voevodsky's theory of homotopy invariant presheaves with transfers relies on an algebro-geometric version of classical homotopy theory, where the r\^ole of the interval is played by the affine line $\A^1$. Reciprocity presheaves with transfers were introduced in \cite{rec} to generalize the former, based on the completely different idea of reciprocity \`a la Rosenlicht-Serre. Conversely, the above theorems say that one may largely understand them in terms of a more sophisticated homotopy theory, based on $\bcube$ rather than $\A^1$. This is a remarkable fact.

\begin{Rm}
In \cite[Conjecture 1 (1)]{rec},
it is conjectured that the Cousin complex attached 
to $F \in \Rec \cap \NST$ is exact.
This is proved for $F \in \RSC \cap \NST$ in \cite[Cor. 3]{shuji}.
Thus, Theorem \ref{C1} (2) permits us to deduce
the full statement of the original conjecture.
\end{Rm}

\subsection*{Corrections} 
In the first version of this paper \cite{ksy2}, we made the following two claims about the functor $\omega_\CI$ from Theorem \ref{T1} (2): it  induces 1) a monoidal structure on $\RSC$ from the one on $\CI$, and 2) an equivalence of categories 
\[\CI\cap \MNST \iso \RSC\cap \NST,\]
where $\MNST\subset \MPST$ is the full subcategory of \emph{`Modulus Nisnevich sheaves with transfers'} (see \S \ref{s1.4}). Both proofs have turned out to be incorrect.
The mistake in 2) originates in a false statement in the initial version of \cite{shuji}, 
which has been removed from its published version. See Remark \ref{rem;counterexample} for a counterexample in characteristic zero. 

\subsection*{Acknowledgments} 
Part of this work was done while the authors stayed at the university of Regensburg supported by the SFB grant ``Higher Invariants". Another part was done in a Research in trio in CIRM, Luminy. Yet another part was done while the third author was visiting IMJ-PRG supported by the Foundation Sciences Ma\-th\'e\-ma\-ti\-ques de Paris. Yet another part was done while the first author was visiting Toh\^oku University and the Tokyo Institute of Technology. We are grateful to the support and hospitality received in all places.

The authors express their deep gratitude 
to Kay R\"ulling for pointing out
the errors mentioned in the previous paragraph.
They are also grateful to the referee for careful reading and useful comments.

\subsection*{Notation and conventions}
Throughout this paper we work over a base field $k$. 
Denote by $\Sch$ the category of separated schemes of finite type over $k$,
and by $\Sm$ the full subcategory of $\Sch$
consisting of all smooth $k$-schemes.

\enlargethispage*{20pt}

\section{Review of basic definitions and results}

\subsection{Modulus pairs}
The following definitions (1) and (2) are taken from \cite[Definitions 1.1.1, 1.3.1]{MwM1}.

\bigskip

\begin{definition}\label{def:mcor}\
\begin{enumerate}
\item
A pair $M=(\ol{X}, X_\infty)$ of 
$\ol{X} \in \Sch$ and 
an effective Cartier divisor $X_\infty$ on $\ol{X}$
is called a \emph{modulus pair}
if $\ol{X} \setminus |X_\infty| \in \Sm$.
It is called {\it proper} if $\ol{X}$ is proper over $k$.
\item
Let 
$M=(\ol{X}, X_\infty), ~N=(\ol{Y}, Y_\infty)$
be two proper modulus pairs and put 
$X=\ol{X} \setminus |X_\infty|, ~Y=\ol{Y} \setminus |Y_\infty|$.
We define $\MCor(M, N)$ to be the subgroup
of $\Cor(X, Y)$ generated by all
elementary correspondences $V \in \Cor(X, Y)$
such that 
the closure $\ol{V}$ of $V$ in $\ol{X} \times \ol{Y}$ satisfies
$\nu^*(\ol{X} \times Y_\infty)\le \nu^*(X_\infty \times \ol{Y})$, where $\nu : \ol{V}^N \to \ol{X}\times \ol{Y}$ is the 
composition of the normalization $ \ol{V}^N \to \ol{V}$ 
and the inclusion $\ol{V} \hookrightarrow \ol{X}\times \ol{Y}$.
We call these correspondences \emph{admissible} (with respect to $(M,N)$).
This defines a category $\MCor$
of proper modulus pairs.
\end{enumerate}
\end{definition}

There is a functor 
\begin{equation}\label{eq;omega}
\omega : \MCor \to \Cor
\end{equation}
defined by $\omega(\ol{X}, X_\infty)=\ol{X} \setminus |X_\infty|$.

The following is basic to his paper.

\begin{lemma}[\protect{\cite[\S 2.1]{MwM3}}]\label{d1.1}
The assignment
\[ (\ol{X}, X_\infty) \otimes (\ol{Y}, Y_\infty)
=(\ol{X} \times \ol{Y}, 
X_\infty \times \ol{Y}+\ol{X} \times Y_\infty).
\]
defines a symmetric monoidal structure on $\MCor$, with unit object $\un=(\Spec k,\emptyset)$. The functor $\omega$ of \eqref{eq;omega} is symmetric monoidal.
\end{lemma}


\subsection{Modulus presheaves with transfers}\label{sect:mpst}

Here is the definition of our main object of study (see {\cite[Definition 2.1.1, Notation 2.1.2]{MwM1}}).

\begin{definition}\label{d2.1}\
\begin{enumerate}
\item We denote by $\MPST$ 
the abelian category of  all additive functors $\MCor^\op \to \Ab$.
\item For $M \in \MCor$, we denote
by $\Z_\tr(M) \in \MPST$ 
the object represented by $M$.
\end{enumerate}
\end{definition}

By \cite[Def. 8.2]{mvw} and \cite[Appendix]{KY} we have the following.

\begin{prop}\label{prop:tensor-mpst}
The category $\MPST$ has a symmetric monoidal structure
that extends the tensor structure of Lemma \ref{d1.1} via the additive Yoneda functor.
It admits an internal Hom such that
\begin{equation}\label{eq;intHomMPST}
 \uHom_{\MPST}(\Z_\tr(M), F)(N)=F(M \otimes N)
\end{equation}
for $M, N \in \MCor$ and $F \in \MPST$.
\end{prop}

\subsection{Relation with $\PST$}\label{s2.3}

The functor $\omega$ of \eqref{eq;omega} induces a functor
$\omega^* : \PST \to \MPST$, $\omega^*(F)=F \circ \omega$.

\begin{proposition}\label{prop:omega-monoidal}\
\begin{enumerate}
\item 
The functor $\omega^*$ is
fully faithful and exact.
\item 
There is a left adjoint 
$\omega_! : \MPST \to \PST$ of $\omega^*$,
which is monoidal and exact.
We have 
\[ \omega_!F(X) \simeq \colim_{M\in \MP(X)} F(M)
\qquad (F \in \MPST, ~X \in \Sm).
\]
\end{enumerate}
\end{proposition}

In (2), 
$\MP(X)$ is the inverse system $\{M=(\ol{X}, X_\infty) \in \MCor\mid X=\ol{X} \setminus |X_\infty|\}$, where transition maps are given by
the diagonal $X \subset X \times X$
whenever it defines a morphism in $\MCor$.

\begin{proof}
See \cite[Prop. 2.2.1 and (2.2.1)]{MwM1}. The monoidality of $\omega_!$ follows that of \eqref{eq;omega}.
\end{proof}

\subsection{Modulus sheaves with transfers}\label{Modulussheaves}\label{s1.4}
In \cite[Lemma-Definition 4.2.1]{MwM2}, we define a full subcategory 
$\MNST\subset \MPST$ 
of ``modulus Nisnevich sheaves with transfers''.
In this paper we need the following:

\begin{prop}\label{p2.1}\
\begin{enumerate}
\item The category $\MNST$ is abelian; the full embedding $i_\Nis:\MNST\inj \MPST$ has an exact left adjoint $a_\Nis$ (``sheafification'').
\item The functors $\omega_!$ and $\omega^*$ of Proposition \ref{prop:omega-monoidal} preserve $\MNST$ and $\NST$; they induce an adjunction $(\omega_\Nis,\omega^\Nis)$ between these two categories, and $\omega_\Nis,\omega^\Nis$ are both exact. Moreover, the pair $(\omega_!,\omega_\Nis)$ 
commutes with the sheafification functors $a_\Nis$ and $a_\Nis^V:\PST\to \NST$ 
\cite[Th. 3.1.4]{voetri}.
\end{enumerate}
\end{prop}

\begin{proof} See \cite[Theorem 4.2.4]{MwM2} for (1) and \cite[Prop. 6.2.1]{MwM2} for (2).
\end{proof}

\section{$\bcube$-invariance and SC-reciprocity}

\subsection{$\bcube$-invariance}

\begin{definition}\label{def:cube-inv} Let $\bcube = (\P^1,\infty)$, and write $p:\bcube\to \un$ for the canonical morphism.
We say 
$F \in \MPST$ is \emph{$\bcube$-invariant}
if the projection map  $1_M\otimes p : M \otimes \bcube \to M$ 
induces an isomorphism
$p^* :F(M) \iso F(M \otimes \bcube)$
for any $M \in \MCor$.
We define $\CI$ to be the full subcategory of $\MPST$
consisting of all objects having
$\bcube$-invariance.
\end{definition}

\begin{lemma}\label{lem;cubeinv} 
The category $\CI$ is closed under taking subobjects, quotients and extensions in $\MPST$. 
\end{lemma}
\begin{proof}
Since the zero section $i_0:\un \to \bcube$ is right inverse to $p$,
$p^*: F(M)\to F(M\otimes\bcube)$ is an isomorphism if and only if
$i_0^*: F(M\otimes\bcube)\to F(M)$ is injective. This implies that $\bcube$-invariance is preserved under taking subobjects.
The remaining assertions then follow by the five lemma.
\end{proof}

Consider the multiplication map
\[ \mu : \A^1 \times \A^1 \to \A^1; \quad (x, y) \mapsto (xy), \]
Let $\Gamma\subset \A^1\times\A^1\times\A^1$ be the graph of $\mu$.

\begin{lemma}[\protect{\cite[Lem. 5.1.1]{MwM3}}]\label{claim;lem0;hcube}
We have $\Gamma \in \MCor(\bcube\otimes \bcube, \bcube)$. In other words, the finite correspondence $\mu$ is admissible.
\end{lemma}

\begin{definition}\label{def;hcubeM}
For $F\in \MPST$, define $\hM F\in \MPST$ by:
\begin{equation}\label{eq;hM}
 \hM F(M) =\Coker\big(F(M\otimes\bcube)\rmapo{i_0^*-i_1^*} F(M)\big)
\;\;(M\in \MCor)\\
\end{equation} 
where $i_\ep^*$ for $\ep=0,1$ is the pullback by the section $i_\ep: \un \to \bcube$ sending $\Spec k$ to $\epsilon\in \A^1(k)$. 
For $M\in \MCor$, we write $\hM M=\hM{\Ztr(M)}$.
\end{definition}

\begin{prop}\label{lem0;hcube}
Let $F\in \MPST$.
\begin{enumerate}
\item
The following conditions are equivalent.
\begin{thlist}
\item
$F\in \CI$.
\item
The natural map $F\to \hM F$ is an isomorphism. 
\item
For any $M\in \MCor$ and $a\in F(M)$, the Yoneda map 
$\ta:\Ztr(M) \to F$ factors through $\hM M$.
\end{thlist}
\item
We have $\hM F\in \CI$ and the induced functor 
\[h_0^{\bcube} :\MPST \to \CI;\;  F\mapsto \hM F\]
gives a left adjoint of the inclusion $i^\bcube:\CI \inj \MPST$. 
\item For any $M\in \MCor$, the morphism $h_0^\bcube(1_M\otimes p):h_0^\bcube(M\otimes \bcube)\to h_0^\bcube(M)$ is an isomorphism.
\end{enumerate}
\end{prop}
\begin{proof} It essentially reproduces the proof of the same facts for $\A^1$-invariant presheaves, by adding modulus. The main point is Lemma \ref{claim;lem0;hcube}.

Assume (i) and take $M\in \MCor$.
The assumption implies that $i_\epsilon^*: F(M\otimes\bcube) \to F(M)$ for $\epsilon=0,1$ are both inverse to $p^*: F(M)\iso F(M\otimes\bcube)$ so that $i_0^*-i_1^*=0$, which implies (ii).

Assume (ii). By \eqref{eq;hM} this implies that for any $M\in \MCor$ we have
\begin{equation}\label{eq1;lem0;hcube}
 i_0^*=i_1^* : F(M\otimes\bcube)\to F(M).
\end{equation}

By Lemma \ref{claim;lem0;hcube}, we have a commutative diagram
\[\xymatrix{
F(M \otimes \bcube) \ar[rr]^{(1_M\otimes i_0)^*} \ar[d]_{(1_M \otimes \mu)^*} && F(M) \ar[d]^{p^*} \\
F(M\otimes \bcube\otimes\bcube) \ar[rr]^{(1_{M\otimes \bcube} \otimes i_0)^*}  &&F(M \otimes \bcube).
}\]
By this diagram and \eqref{eq1;lem0;hcube}, we get
\begin{multline*}
p^* (1_M\otimes i_0)^* 
= (1_{M\otimes \bcube} \otimes i_0)^* \circ (1_M \otimes \mu)^* \\
= (1_{M\otimes \bcube} \otimes i_1)^* \circ (1_M \otimes \mu)^* 
= 1_{M \otimes \bcube}^*.
\end{multline*}
This proves the surjectivity of $p^*$, hence (i) holds. Thus (i) $\iff$ (ii).

By the definition of $\hM F$, for any $M\in \MCor$, the map
\[ \hM F(M\otimes\bcube) \rmapo{i_0^*-i_1^*} \hM F(M) \]
is the zero map so that 
$\hM F(M)\simeq \hM{\hM F}(M)$.
Hence the first assertion of (2) follows from the implication (ii) $\Rightarrow$ (i). 

Any $\ta:\Z_\tr(M)\to F$ induces a morphism $\hM M \to \hM F$ which commutes with the natural transformation of (ii). Hence (iii) follows from (ii).

If (iii) holds, $F$ is a quotient of a direct sum of $\hM M$'s for $M\in \MCor$. Hence (i) holds by the first assertion of (2) and Lemma \ref{lem;cubeinv}. This completes the proof of (1).

To show the second assertion of (2), note that (1) implies that
for $F\in \CI$ and $M\in \MCor$, the natural map $\Ztr(M)\to \hM M$ induces an isomorphism
\begin{equation}\label{eq;hcube}
\Hom_{\MPST}(\hM M,F) \simeq \Hom_{\MPST}(\Ztr(M),F).
\end{equation}
For $G\in \MPST$, take a resolution of $G$ of the form
\[ P_1 \to P_0 \to G\to 0 \;\;\text{ in $\MPST$,}\]
where $P_1,P_0$ are direct sums of representable objects.
By its definition \eqref{eq;hM}, the endofunctor $h_0^{\bcube}$ of $\MPST$ is right exact so that the above sequence induce an exact sequence
\[ \hM{P_1} \to \hM{P_0} \to \hM{G}\to 0.\]
Moreover $h_0^{\bcube}$ commutes with direct sums again by \eqref{eq;hM}.
In view of \eqref{eq;hcube}, we conclude that the natural map $G\to \hM G$ induces an isomorphism
\[ \Hom_{\MPST}(\hM G,F) \simeq \Hom_{\MPST}(G,F).\]
which implies the desired claim. 

It remains to prove (3). Since $h_0^\bcube(1_M\otimes i_0)$ is a right inverse of $h_0^\bcube(1_M\otimes p)$, it suffices to show that it is also a left inverse. Let $N\in \MCor$. 
For $\phi\in \MCor(N,M\otimes \bcube)$, define
\[\tilde \phi = (1_M\otimes \mu)\circ (\phi\otimes 1_\bcube)\in \MCor(N\otimes \bcube, M\otimes \bcube).\]

Using the identities $\mu\circ (1_\bcube\otimes i_0)=i_0p$ and $\mu\circ (1_\bcube\otimes i_1)=1_\bcube$, we get
\[(i_0^*-i_1^*)\tilde \phi :=\tilde\phi\circ (1_N\otimes i_0)-\tilde\phi \circ (1_N\otimes i_1) = (1_M\otimes  i_0p)\circ \phi -\phi  \]
which shows that $\bar \phi=h_0^\bcube(1_M\otimes i_0)h_0^\bcube(1_M\otimes p)\bar \phi$, where $\bar \phi$ is the image of $\phi$ in $h_0^\bcube(M\otimes \bcube)(N)$. 
\end{proof}

\begin{definition}\label{def;hMMcube}
For $F\in \MPST$, define $\hMM F\in \MPST$ by
\begin{equation}\label{eq;hMM}
\hMM F(M) = \Hom_{\MPST}(\hM {M},F)\;\; (M\in \MCor).
\end{equation}
\end{definition}

\begin{lemma}\label{lem;hMM}
For $F\in \MPST$, $\hMM F$ is the maximal $\bcube$-invariant subobject of $F$. 
The induced functor 
\[h^0_{\bcube} :\MPST \to \CI;\;  F\mapsto \hMM F\]
gives a right adjoint of the inclusion $i^\bcube:\CI \inj \MPST$. 
\end{lemma}
\begin{proof}
The fact that $\hMM F$ is a subobject of $F$ follows from \eqref{eq;hMM} and the fact that $\hM M$ is a quotient of $\Ztr(M)$.
The fact that $\hMM F\in \CI$ follows from 
Proposition \ref{lem0;hcube} (3). Now let $G\subset F$ be a subobject which is in $\CI$.
For $a\in F(M)$ with $M\in \MCor$, let $\ta:\Ztr(M)\to F$ be the corresponding map in $\MPST$ . If $a\in G(M)$, $\ta$ factors through $G$ and hence
factors through $\hM M$ by Proposition \ref{lem0;hcube} (1). Hence $a\in \hMM F(M)$ by \eqref{eq;hMM}. This proves $G\subset \hMM F$, which completes the proof of the first assertion. The second assertion follows easily from the first and Lemma \ref{lem;cubeinv}.
\end{proof}

\enlargethispage*{30pt}

\begin{thm}\label{rem:monoidal-ci}\
The category $\CI$ is a Serre subcategory of $\MPST$,  and is Grothendieck.
The inclusion $i^{\bcube} :\CI \inj \MPST$ has
\begin{thlist}
\item a left adjoint given by 
$F \mapsto h_0^{\bcube}(F)$;
\item a right adjoint given by $F\mapsto h^0_{\bcube}(F)$, where 
\[h^0_{\bcube} F(M) = \Hom(h_0^{\bcube}(M),F)\quad (M\in \MCor).\]
\end{thlist}
The unit (resp. counit) morphism $F\to h_0^\bcube(F)$ (resp. $h^0_\bcube(F)\to F$) is epi (resp. mono).
\end{thm}

 \begin{proof} Everything follows from what we have proven so far, except for the Grothendieckness of $\CI$, that we want to deduce from that of $\MPST$. We cannot quite apply \cite[Th. A.10.1 d)]{MwM1}, because $h_0^\bcube$ is not exact: it provides generators and infinite direct sums, but not their exactness. But the latter holds because $i^\bcube$ reflects infinite direct sums (if $(F_\alpha)$ is a family of objects in $\CI$, their direct sum in $\MPST$ belongs to $\CI$).
\end{proof}

\begin{prop}\label{prop;CIotimes}\
\begin{enumerate}
\item One has $\uHom_\MPST(G,H)\in \CI$  for any $G\in \MPST$ and any $H\in \CI$.
\item Via $h_0^\bcube$, the symmetric monoidal structure on $\MPST$ from Proposition \ref{prop:tensor-mpst} induces a symmetric monoidal structure on $\CI$ by the formula
\[ F \otimes_{\CI} G =
h_0^{\bcube}(F \otimes_{\MPST} G) \qquad (F, G \in \CI).
\]
This tensor product commutes with all representable colimits.
\end{enumerate}
\end{prop}

\begin{proof} (1) follows from  \eqref{eq;intHomMPST} and the isomorphisms (easily checked by means of Yoneda's lemma):
\begin{multline*}
\uHom_\MPST(\bcube,\uHom_\MPST(G,H))\simeq 
\uHom_\MPST(\bcube \otimes_\MPST G,H)\\
\simeq \uHom_\MPST(G,\uHom_\MPST(\bcube,H))\overset{\sim}{\gets}
\uHom_\MPST(G,H)
\end{multline*}
where the last isomorphism, induced by $p$, follows from the assumption $H\in \CI$.

(2) Since $h_0^\bcube$ is a localisation by the full faithfulness of $i^\bcube$, we have to show the following.

\begin{claim}
If $f\in \Hom_{\MPST}(G_1,G_2)$ is such that $\hM {f}$ is an isomorphism, 
then $g=\hM{f\otimes_\MPST 1_{G'}}$ is an isomorphism for any $G'\in \MPST$. 
\end{claim}

By (co)Yoneda, it suffices to show that
\[g^*:\Hom_{\MPST}(\hM {G_2\otimes_\MPST G'},H)\to 
\Hom_{\MPST}(\hM{G_1\otimes_\MPST G'},H) \]
is an isomorphism for any $H\in \CI$. By adjunction, it suffices to show that 
\[f^*:\Hom_{\MPST}(G_2,\uHom_\MPST(G',H))\iso 
\Hom_{\MPST}(G_1,\uHom_\MPST(G',H))\]
which follows from (1).

Finally, the commutation of $\otimes_\CI$ with colimits follows from that of $\otimes_\MPST$ and $h_0^\bcube$.
\end{proof}

\subsection{$SC$-reciprocity}\label{sect:rel-rec}

\begin{definition}\label{d4.3}
For $F \in \MPST$, we define 
\[ \hMw{F}:=\omega_! \hM {F} \in \PST.
\]
For $M \in \MCor$, we put $\hMw {M}:= \hMw {\Z_{\tr}(M)}$.
\end{definition}

\begin{lemma}\label{lem;hMw}
Let $M=(\Mb,\Minf) \in \MCor$ with $M^\o=M-|\Minf|$.
For $S\in \Sm$, we have
\begin{equation*}\label{eq;d4.3}
\hMw {M}(S)=\Coker(i_0^{*} - i_1^{*}:\ulMCor(\bcube\otimes S, M) \to \Cor(S,M^\circ)\big),
\end{equation*}
where $\ulMCor(\bcube\otimes S, M)$ is the subgroup of 
$\Cor(\A^1\times S, M^\circ)$ generated by all elementary correspondences $Z$ such that 
\[
\phi_Z^*(\P^1\times S \times \Minf ) \leq
\phi_Z^*(\infty \times S \times \Mb),
\]
where $\phi_Z : \ol{Z}^N \to \ol{Z} \hookrightarrow \P^1 \times S \times \Mb$
denotes the normalization of the closure $\ol{Z}$ of $Z$ in $\P^1 \times S \times \Mb$. 
\end{lemma}
\begin{proof}
This follows from Definition \ref{def;hcubeM} and 
 Proposition \ref{prop:omega-monoidal} (2). 
\end{proof}

\begin{remark}\label{rem:sushom-chow}
Lemma \ref{lem;hMw} implies isomorphisms
\[ \hMw {M}(\Spec k) \simeq  H_0^S(\Mb, \Minf) \simeq \CH_0(M),\]
where $H_n^S(\Mb, \Minf)$ for $n\in \Z$ is the {\it Suslin homology} 
considered in \cite[Definition 3.1]{ry} and 
$\CH_0(M)=\CH_0(\Mb|\Minf)$ is the Chow group with modulus considered in \cite{KSa}.
\end{remark}

\begin{definition}\label{d2.3-2}\
\begin{enumerate}
\item
Let $F \in \PST$,
$X \in \Sm$ and $a \in F(X)=\Hom_\PST(\Z_\tr(X), F)$.
We say $M=(\ol{X}, X_\infty) \in \MCor$ is an \emph{SC-modulus for $a$} if
$X=\ol{X} \setminus |X_\infty|$ and
$a : \Z_\tr(X) \to F$
factors through $\Z_\tr(X) \twoheadrightarrow h_0(M)$.
(SC stands for ``Suslin complex''.)
\item
We say $F \in \PST$ \emph{has SC-reciprocity} if,
for any $X \in \Sm$, any $a \in F(X)$ has an SC-modulus $M \in \MP(X)$.
\item
We define $\SCRec$ to be the
full subcategory of $\PST$
consisting of all objects having SC-reciprocity.
\end{enumerate}
\end{definition}

\begin{remark}\label{rem:closed}
The category $\SCRec$ is closed under
subobjects and quotient objects in $\PST$.
This is obvious from the definition. In particular, $\RSC$ is abelian and the inclusion functor $i^\natural:\RSC \to \PST$ is exact.
\end{remark}

Recall that $i^{\bcube} :\CI \inj \MPST$ denotes the inclusion.

\begin{prop}\label{pRR2} The functor $\rho:=\omega_! i^\bcube h_\bcube^0 \omega^*$ sends $\PST$ into $\RSC$, and is right adjoint to the inclusion $i^\natural:\RSC\inj \PST$.
\end{prop}

\begin{proof} For $F\in \PST$ and $X\in \Sm$, we have by successive adjunctions and by Proposition \ref{prop:omega-monoidal} (2):
\begin{multline}\label{eqRR1}
\rho F(X) = \colim_{M\in \MP(X)} i^\bcube h_\bcube^0 \omega^*F(M)\\
 =\colim_{M\in \MP(X)} \MPST(h^\bcube_0(M), \omega^*F) = \colim_{M\in \MP(X)} \PST(h_0(M),F)
\end{multline}
which realises $\rho F$ as the largest subobject of $F$ which is in $\RSC$.
\end{proof}

\begin{cor}\label{lem:flat-star-sh5}
Let $F \in \PST$.  The counit map
\begin{equation}\label{eq:rec-adj5}
i^\natural \rho F \to F
\end{equation}
of the adjunction in Proposition \ref{pRR2} agrees with the counit map of the adjunction $(\omega_!i^\bcube,h^0_\bcube\omega^*)$. Moreover, 
$F\in \SCRec$
if and only if \eqref{eq:rec-adj5} is an isomorphism.
\end{cor}

\begin{proof}
The first statement simply restates the computation in \eqref{eqRR1}. The second one follows from Proposition \ref{pRR2} and the definition of SC-reciprocity.
\end{proof}

\enlargethispage*{30pt}

\begin{qn}
Is $\RSC$ closed under extensions in $\PST$?
\end{qn}

This question appears to be very difficult (compare \cite[Question 1]{rec}). We can only offer a trivial reduction:

\begin{prop}\label{thm;RSCext}
For  $F\in \PST$, the following are equivalent:
\begin{thlist}
\item There exist  $G, H\in \RSC$ and an exact sequence in $\PST$
\begin{equation}\label{eq1;RSCext}
0\to G \to F \to H \to 0.
\end{equation}
\item The cokernel of \eqref{eq:rec-adj5} is in $\RSC$.
\end{thlist}
\end{prop}

\begin{proof} (ii) $\Rightarrow$ (i) is obvious. Conversely, assume (i). Applying the left exact functor $\rho$, we get an exact sequence in $\RSC$ (defining $C$)
\begin{equation}\label{eq2;RSCext**}
 0\to G \to \rho F \to  H\to  C  \to 0.
\end{equation}
The counit map (2.6) sends \eqref{eq2;RSCext**} to \eqref{eq1;RSCext}. A diagram chase then gives us the exact sequence in $\PST$
\[0\to \rho F\to F \to  C\to 0\]
which concludes the proof.
\end{proof}

\subsection{Relations between $\CI$, $\HI$ and $\SCRec$}\label{s2.2}

Recall that  $\HI\subset \PST$ is the full subcategory of $\A^1$-invariant presheaves with transfeers.
 
\begin{lemma}\label{lem:2cube-inv}\
For $H\in \PST$, $H\in \HI$ if and only if $\omega^*H\in \CI$.
\end{lemma}
\begin{proof}
This follows from the fact that 
for $M=(\Mb,\Minf)\in \MCor$ with $M^\circ=\Mb-\Minf$, we have
\[\omega^*F(M)=F(M^\circ) \qaq
\omega^*F(M\otimes\bcube)=F(M^\circ\times\A^1).\]
 \end{proof}

\begin{prop}\label{p4.1} 
The composite $\MPST \rmapo{\omega_!} \PST \rmapo{h_0^{\A^1}} \HI$ factors
through $\MPST \rmapo{h_0^\bcube} \CI$, inducing a functor $\omega_h:\CI\to \HI$.
This functor is right exact and monoidal for the $\otimes$-structures given on $\CI$ by Theorem \ref{rem:monoidal-ci} (2), and analogously on $\HI$. It has an exact right adjoint $\omega^h$, given by the restriction of $\omega^*$ to $\HI$. \end{prop}

\begin{proof} The first claim and the monoidality of $\omega_h$ follow from that of $\omega_!$, as $\omega_! \Z_\tr(\bcube)=\Z_\tr(\A^1)$. The existence, characterisation and exactness of $\omega^h$ follows from Lemma \ref{lem:2cube-inv}, and the right exactness of $\omega_h$ then follows.
\end{proof}

\begin{thm}\label{cor:cube-inv-rec2}
If $F \in \CI$, then $\omega_!F \in \RSC$.
\end{thm}

\begin{proof} We have a commutative diagram, for any $F\in \MPST$:
\begin{equation}\label{eq4.3}
\begin{CD}
\omega_!i^\bcube h^0_\bcube F@>\omega_!i^\bcube h^0_\bcube\eta_F>(a)>  \omega_! i^\bcube h^0_\bcube\omega^*\omega_! F \\
@V \omega_!\epsilon'_FV(c)V @V\omega_!\epsilon'_{\omega^* \omega_!F}V(d)V\\
\omega_! F @>\omega_! \eta_F>(b)> \omega_!\omega^* \omega_! F @>\epsilon_{\omega_!F}>(e)> \omega_!F.
\end{CD}
\end{equation}

Here, $\eta$ and $\epsilon$ are the unit and counit of the adjunction $(\omega_!,\omega^*)$, while $\epsilon'$ is the counit of the adjunction $(i^\bcube,h^0_\bcube)$. We have $(e) \circ (b) = 1_{\omega_!F}$ by the adjunction identities; since $\omega^*$ is fully faithful, (e) is an isomorphism hence so is (b). This shows that (c) factors through (a). On the other hand, $\epsilon'$ is mono by Theorem \ref{rem:monoidal-ci}, hence so are (c) and (d) since $\omega_!$ is exact. Finally, the diagram boils down to two successive monomorphisms
\begin{equation}\label{eq4.1}
\omega_!i^\bcube h^0_\bcube F\Inj i^\natural\rho \omega_! F \Inj \omega_! F 
\end{equation}
with composition $\omega_! \epsilon'_F$. Therefore, $F\in \CI$ $\Rightarrow$ $\omega_!F\in \RSC$.
\end{proof}

\begin{cor}\label{c4.2} We have $\HI\subset \RSC$.
\end{cor}

\begin{proof} Let $F\in \HI$. By Lemma \ref{lem:2cube-inv}, $\omega^*F\in \CI$, hence
\[F\simeq \omega_!\omega^*F\in \RSC\]
by Theorem \ref{cor:cube-inv-rec2}. (See \cite[Lemma 1.22]{shuji} for a simpler proof.)
\end{proof}

\begin{corollary}\label{cor:sus-hom-screc}
For any $F \in \MPST$, $\hMw {F}\in \RSC$.
\end{corollary}
\begin{proof}
This follows from Proposition \ref{lem0;hcube} and
Theorem \ref{cor:cube-inv-rec2}.
\end{proof}

\begin{cor} The inclusion functor $i^\natural:\RSC\inj \PST$ has a pro-left adjoint $\ell$.
\end{cor}

\begin{proof} It suffices to show that $\ell$ is defined on the generators $\Z_\tr(X)$. Since $h_0(M)\in \RSC$ for any $M\in \MP(X)$ by Corollary \ref{cor:sus-hom-screc}, we have $\ell \Z_\tr(X) = ``\lim"_{M\in \MP(X)} h_0(M)$.
\end{proof}

\begin{proposition}\label{prop:omega-ci}
There exist unique functors
$\omega_{\CI}$ and $\omega^{\CI}$
that make the two diagrams
\[
\xymatrix{
\CI
\ar@^{{(}->}[r]^{i^{\bcube}}
\ar[d]_{\omega_{\CI}}
&
\MPST 
\ar[d]_{\omega_!}
&
\CI
&
\MPST 
\ar[l]_{h^0_{\bcube}}
\\
\SCRec
\ar@^{{(}->}[r]_{i^{\natural}}
&
\PST
&
\SCRec
\ar@^{{(}->}[r]_{i^{\natural}}
\ar[u]^{\omega^{\CI}}
&
\PST
\ar[u]_{\omega^*}
}
\]
commutative, where 
$i^{\natural}$ is the inclusion.
Moreover,
$\omega^{\CI}$ is right adjoint to $\omega_{\CI}$. The counit map
$\epsilon:\omega_{\CI} \omega^{\CI} \Rightarrow \id_{\SCRec}$
is an isomorphism, $\omega_{\CI}$ is a localisation (in particular, is essentially surjective) and
$\omega^{\CI}$ is fully faithful. Finally, $\omega_{\CI}$ is exact and $\omega^{\CI}$ is left exact.
\end{proposition}

\begin{proof}
The existence of $\omega_{\CI}$ is the contents of Theorem \ref{cor:cube-inv-rec2}, and $\omega^{\CI}$ is defined by the commutativity of the diagram. 
For the second assertion, let $F\in \CI$ and $G\in \RSC$. Using two successive adjunctions, we compute:
\begin{multline*}
\CI(F,\omega^{\CI} G) = \CI(F,h^0_\bcube \omega^*i^\natural G)\simeq \PST(\omega_! i^\bcube F, i^\natural G)\\ =\PST(i^\natural \omega_{\CI}  F, i^\natural G)\simeq \RSC(\omega_{\CI}  F, G)
\end{multline*}
where the last isomorphism uses the (tautological) full faithfulness of $i^{\natural}$. So the adjunction $(\omega_{\CI}, \omega^{\CI})$ is obtained by ``cancelling'' $i^\natural$ from the adjunction $(\omega_! i^\bcube,h^0_\bcube\omega^*)$, after applying Theorem \ref{cor:cube-inv-rec2}. Therefore the third assertion follows from Corollary \ref{lem:flat-star-sh5}, 
and the next two are standard consequences \cite[Lemma A.3.1]{MwM1}. The exactness of $\omega_{\CI}$ follows from the exactness of $i^\bcube$ and $\omega_!$ (as well as the full faithfulness of $i^\natural$), and $\omega^{\CI}$ is left exact as a right adjoint.
\end{proof}

\begin{cor}\label{c2.1} The category $\RSC$ is Grothendieck.
\end{cor}

\begin{proof} This follows from the same fact for $\CI$ (Theorem \ref{rem:monoidal-ci}), the adjunction $(\omega_\CI,\omega^\CI)$ and  \cite[Th. A.10.1 d)]{MwM1}.
\end{proof}

\begin{prop}\label{p4.2}
Let 
\[ h_0^{\rec}:\RSC\to \HI\]
be the restriction of $h_0^{\A^1}:\PST\to \HI$ from \eqref{h0A}.
Then $h_0^{\rec}$ is a left adjoint of the inclusion $\HI\hookrightarrow \RSC$ from Corollary \ref{c4.2}.
We have a natural isomorphism $\omega_h\simeq h_0^{\rec} \omega_{\CI}$ (see Proposition \ref{p4.1} for $\omega_h$). 
\end{prop}

\begin{proof} The first claim follows immediately from the fact that 
$h_0^{\A^1}$ is a left adjoint to the inclusion $\HI\hookrightarrow \PST$.
To show the second, we apply the natural isomorphism $\omega_h h_0^\bcube G\simeq h_0^{\A^1} \omega_! G$ from Proposition \ref{p4.1} to $G=i^\bcube F$ for $F\in \CI$ to get a natural isomorphism 
\[\omega_hF \simeq \omega_h h_0^\bcube i^\bcube F \simeq h_0^{\A^1}\omega_! i^\bcube F \simeq h_0^{\A^1}i^\natural \omega_{\CI}F\simeq h_0^{\rec} \omega_{\CI}F\]
as requested. 
\end{proof}

\subsection{Sheaves in $\RSC$}\label{s4.7a} 

Let $\NST\subset \PST$ be the full subcategory of Nisnevich sheaves with transfers
\cite[Th. 3.1.4]{voetri}. Recall that the objects of $\NST$ are those $F\in \PST$ whose restriction $F_X$ to $X_\Nis$ is a sheaf for any $X \in \Sm$,
where $X_\Nis$ denotes the small Nisnevich site of $X$. 
By \cite[Th. 3.1.4]{voetri} the inclusion $i^V_\Nis:\NST\to \PST$ has an exact left adjoint $\aVNis$ such that for any $F\in \PST$ and $X\in \Sm$, 
$(\aVNis F)_X$ is the Nisnevich sheafication of $F_X$ as a presheaf on $X_\Nis$. 
\def\RSCNis{\RSC_{\Nis}}
Let $\RSCNis=\RSC\cap \NST$ and 
$\CI_\Nis=\CI\cap \MNST$ (see \S \ref{s1.4} for $\MNST$).
We admit the following theorem.

\begin{thm}\label{t3.1}
Assume $k$ is perfect. Write
\[\CI^{sp} = \{F\in \CI\mid \text{ the unit map } F\to \omega^\CI\omega_\CI F \text{is injective.}\}
\]
\begin{enumerate}
\item \cite[Th. 0.1 and 0.4]{shuji} 
One has $\aVNis(\RSC)= \RSC_\Nis$
and $a_\Nis \CI^{sp} \subset \CI_\Nis$.  (See Proposition \ref{p2.1} (1) for $a_\Nis$.)
\item \cite[Cor. 4.16]{rs2}. One has $\omega^\CI(\RSC_\Nis) \subset \CI_\Nis$. 
\end{enumerate}
\end{thm}

\begin{cor} The category $\RSC_\Nis$ is Grothendieck.
\end{cor}

\begin{proof} Since $a_\Nis^V$ is exact, so is its restriction to $\RSC$. The corollary now follows from Corollary \ref{c2.1} and (again)   \cite[Th. A.10.1 d)]{MwM1}.
\end{proof}

\begin{thm}\label{cRR1} \
Assume $k$ is perfect.
\begin{enumerate}
\item The functor $\rho$ of Proposition \ref{pRR2} sends $\NST$ into $\RSC_\Nis$. It yields a right adjoint $\rho_\Nis$ to the inclusion $i^\natural_\Nis:\RSC_\Nis\inj \NST$.
\item 
The functor $\omega_{\CI}$ of Proposition \ref{prop:omega-ci} sends 
$\CI_\Nis$ to $\RSC_\Nis$.
The induced functor $\omega_\CI^\Nis :\CI_\Nis \to \RSC_\Nis$  is left adjoint to the fully faithful functor $\omega^\CI_\Nis:\RSC_\Nis \to \CI_\Nis$ given by Theorem \ref{t3.1} (2).   
Moreover, there is a natural ismorphism
\begin{equation}\label{eq3.6}
a_\Nis^V \omega_\CI F\simeq \omega_\CI^\Nis a_\Nis F
\end{equation}
for any $F\in  \CI^{sp}$.
\end{enumerate}
\end{thm}

\begin{proof} 
Let $F\in \NST$. Considering $F$ as an object of $\PST$, we may view 
$\rho F$ as the largest subobject of $F$ which belongs to $\RSC$ (see Proposition \ref{pRR2}). Applying the left exact functor $\aVNis$ to this inclusion, we get a sequence
\[\rho F\to \aVNis \rho F\to \aVNis F=F,\]
where the second map is a monomorphism. But the middle term is in $\RSC$ by Theorem \ref{t3.1} (1). Hence the first map must be an isomorphism, which implies the first claim of (1). The last claim now follows easily from the adjunction in Proposition \ref{pRR2}.

The first assertion of (2) is obvious 
since $\omega_!$ preserves Nisnevich sheaves by Proposition \ref{p2.1}.  The second one then follows easily from Proposition \ref{prop:omega-ci}. Given the natural isomorphism $\omega^\CI i_\Nis^V\simeq i_\Nis \omega^\CI_\Nis$,
this implies the last assertion by taking left adjoints.
\end{proof}

\begin{remark}\label{rem;counterexample}
The functor $\omega_{\CI}^\Nis$ is not conservative.
Assume $\ch(k)=0$.  Let $F\in \CI$ be the image of the unit map
\[ h_0^{\bcube}(\P^1,2\infty) \to \omega^\CI\omega_\CI h_0^{\bcube}(\P^1,2\infty).\]
Then $F\in \CI^{sp}$, hence $a_\Nis F\in \CI$ by Theorem \ref{t3.1} (1). We claim that the unit map
$\iota : a_\Nis F \to \omega^{\CI}_\Nis\omega_{\CI}^\Nis a_\Nis F$
 is not surjective. To see this, first note that 
by the exactness of $\omega_\CI$, we have
\[ h_0(\P^1,2\infty) 
\twoheadrightarrow \omega_\CI F
\hookrightarrow
\omega_\CI \omega^\CI h_0(\P^1,2\infty)=h_0(\P^1,2\infty),
\]
and hence $\omega_\CI F \cong h_0(\P^1,2\infty)$.
Then by \eqref{eq3.6}  and
\cite[Thm. 1.1]{ry}, we have 
isomorphisms 
\[
\omega_{\CI}^\Nis a_\Nis F
\simeq 
a_\Nis^V \omega_\CI F
\simeq \underline{\Pic}(\P^1,2\infty)\simeq \Z\oplus \G_a.
\]
Take $(X,D)\in \MCor$ such that $X,D\in \Sm$, with $X$ connected.
Then it follows from \cite[Th. 6.4]{rs2} that
\[\omega^{\CI}_\Nis\omega_{\CI}^\Nis a_\Nis F (X,2mD) \simeq
\Z \oplus H^0(X,\sO_X((2m-1)D))\]
 for any integer $m>0$. On the other hand, one can show 
\[a_\Nis F (X,2m D) \simeq \Z \oplus H^0(X,\sO_X(m D)).\]
This implies that $G:=\Coker(\iota)\in \CI_\Nis$ is non-zero but $\omega_{\CI}^\Nis(G)=0$.
\end{remark}

\section{Relation with \cite{rec}}

\subsection{Review of reciprocity presheaves with transfers}\label{ssec:ref} In \cite[Definition 2.1.3]{rec},
we defined a full subcategory $\Rec$ of $\PST$,
which we now recall.

Let $(\Xb, Y) \in \MCor$ and suppose that $X=\Xb \setminus |Y|$
is quasi-affine.
For $S\in \Sm$, 
let $\sC_{(\Xb, Y)}(S)$ be the  class of all finite morphisms 
$\phi : \ol{C} \to \Xb \times S$
satisfying the following conditions:
\begin{itemize}
\item 
$\ol{C} \in \Sch$ is integral and normal.
\item
There is a generic point $\eta$ of $S$ such that
$\dim \ol{C} \times_S \eta=1$.
\item
The image of $\gamma_\phi := \mathrm{pr} \circ \phi$ is not contained in $|Y|$,
where $\mathrm{pr} : \ol{X} \times S \to \ol{X}$ is the projection map.
\end{itemize}
For an effective Cartier divisor $D$ on $\ol{C}$,
we set
\begin{equation}\label{eq.GCD}
G(\ol{C},D):= \bigcap_{x\in D} \Ker\big(\sO_{\ol{C},x}^\times 
\to \sO_{D,x}^\times\big).
\end{equation}
We then define
\[  \Phi(\Xb, Y)(S) = 
\bigoplus_{(\phi : \ol{C} \to \Xb \times S) \in \sC_{(\Xb, Y)}(S)} 
G(\ol{C}, \gamma_{\phi}^* Y) 
\]
It is proved in \cite[Proposition 2.2.2]{rec} that
$\Phi(\Xb, Y)$ defines a presheaf with transfers.
It is also shown there that
one has $\phi_*(\div_{\ol{C}}(f)) \in \Cor(S, X)$
for any $(\phi : \ol{C} \to \Xb \times S) \in \sC_{(\Xb, Y)}(S)$
and $f \in G(\ol{C}, \gamma_{\phi}^* Y)$,
yielding a map
$\tau : \Phi(\Xb, Y) \to \Z_\tr(X)$ in $\PST$.
We define
\[ h(M):= \Coker(\tau : \Phi(\Xb, Y) \to \Z_\tr(X)) \in \PST. \]

\begin{definition}[{\cite[{Definition 2.1.2, Remark 2.1.6}]{rec}}]
\label{def:old-rec}
We say $F \in \PST$ has reciprocity
if for any quasi-affine $X \in \Sm$ and 
$a \in F(X)=\Hom_\PST(\Z_\tr(X), F)$,
there is an $M=(\ol{X}, X_\infty) \in \MCor$ 
such that
$X=\ol{X} \setminus |X_\infty|$ and
$a : \Z_\tr(X) \to F$
factors through $\Z_\tr(X) \twoheadrightarrow h(M)$.
We define $\Rec$ to be the
full subcategory of $\PST$
consisting of all objects having reciprocity.
\end{definition}

\subsection{Statement of the result and consequences}

\begin{thm}\label{thm:h-h0} Let $M=(\ol{X}, Y) \in \MCor$ be such that
$X:=\ol{X} \setminus |Y|$ is quasi-affine. Then $h_0(M)=h(M)$.
Hence we have $\SCRec \subset \Rec$.
\end{thm}

The proof of 
Theorem \ref{thm:h-h0}  will occupy \S\S \ref{s5.1}
and \ref{subsect:screc-rec}. 
We first deduce some consequences.

\begin{cor}
For any $F \in \RSC$, 
we have $F_\Zar \simeq F_\Nis$,
where $F_\Zar$ (resp. $F_\Nis$)
is the Zariski (resp. Nisnevich) sheafification 
of $F$.
\end{cor}
\begin{proof}
Combine Theorem  \ref{thm:h-h0}
and \cite[Theorem 7]{rec}.
\end{proof}

The next result depends on Theorem \ref{t3.1} (1).

\begin{cor}\label{prop:sc-rec-rec}\
Assume $k$ is perfect.
Then we have $\SCRec_\Nis = \Rec_\Nis$.
\end{cor}

\begin{proof} The inclusion follows immediately from Theorem \ref{thm:h-h0}. To prove the equality,  let  $F\in \Rec_\Nis$. By \eqref{eqRR1} and Theorem \ref{thm:h-h0}, the map $i^\natural\rho F \to F$ of \eqref{eq:rec-adj5} is an isomorphism when evaluated at $X$ if $X$ is quasi-affine. By Theorem \ref{cRR1} (1), this extends to any $X\in \Sm$ by using a quasi-affine Zariski cover. Thus $F\in \RSC_\Nis$.
\end{proof}

\begin{rk} 
Here is an example of an object $F\in \Rec\setminus \RSC$. Define $F$ as
\[\Coker\left(\bigoplus_{(X,a)} \Z_\tr(X)\to \Z_\tr(\P^1)\right)\]
where $X$ runs through all smooth quasi-affine $k$-schemes and $a$ runs through all elements of $\Cor(X,\P^1)$. By construction, $F(X)=0$ for any smooth quasi-affine $X$, hence $F\in \Rec$. On the other hand, we claim that the image $\eta\in F(\P^1)$ of the identity map $1_{\P^1}\in \Z_\tr(\P^1)(\P^1)$ does not have an SC modulus.  Since $\P^1$ is proper, this amounts to say that the composition 
\[\uHom(\Z_\tr(\A^1),\Z_\tr(\P^1))\by{i_0^*-i_1^*} \Z_\tr(\P^1)\by{\eta} F\]
is nonzero. 
The quasi-affineness of the $X$'s yields that 
for any proper $Y \in \Sm$ the image of
$\bigoplus_{(X,a)} \Z_\tr(X)(Y)\to \Z_\tr(\P^1)(Y)=\Cor(Y, \P^1)$
is generated by cycles of the form
$Y \times x$ where $x$ ranges over closed points of $\P^1$.
In particular, if we take $Y=\P^1$
we find that $F(\P^1)$ is not finitely generated.
On the other hand, 
\cite[Th. 3.3.1]{birat-tri} shows
\begin{multline*}
\Coker(\uHom(\Z_\tr(\A^1),\Z_\tr(\P^1))(\P^1)\by{i_0^*-i_1^*} \Z_\tr(\P^1)(\P^1))\\
\simeq \Pic(\P^1 \times \P^1)
\simeq \Z \times \Z.
\end{multline*} 
Hence $\eta$ cannot vanish at $\P^1$.
\end{rk}

\begin{corollary}\label{c3.2}
Assume $k$ is perfect.
\begin{enumerate}
\item 
A presheaf with transfers represented by a
smooth commutative algebraic group has SC-reciprocity.
\item
The presheaf with transfers $H^0(-, \Omega_-^i)$
has SC-reciprocity for any $i \geq 0$.
The same is true for 
the presheaf with transfers $H^0(-, \Omega_{-/k}^i)$.
\item
Suppose that $k$ is of positive characteristic.
Then the presheaf with transfers $H^0(-, W_n \Omega_-^i)$
has SC-reciprocity for any $i \geq 0$ and $n \geq 1$.
\end{enumerate}
\end{corollary}
\begin{proof}
Combine Corollary \ref{prop:sc-rec-rec}
and \cite[Theorems 4, 5]{rec}.
\end{proof}

The next corollary uses the work of Binda et al \cite{many}: we suppose $k$ is of characteristic $p>0$ and
we use the notation $[1/p]$ to designate categories constructed out of sheaves of $\Z[1/p]$-modules: they are full subcategories of those considered in this paper. 

\begin{cor}\label{c5.2} Assume 
that $\car k=p>0$. Then the functor $h_0^{\rec}$ from Proposition \ref{p4.2} induces an equivalence of categories 
\[\RSC_\Nis[1/p]\iso \HI_\Nis[1/p].\] 
\end{cor}

\begin{proof} By Proposition \ref{p4.2}, it suffices to show that $F\iso h_0^\rec(F)$ for any $F\in \RSC_\Nis[1/p]$. If $F\in \Rec_\Nis[1/p]$, this follows from \cite[Th. 3.5 (2)]{many}, hence the claim when $k$ is perfect by  Corollary  \ref{prop:sc-rec-rec}. The general case reduces to this one by \cite[Prop. 4.5]{adjoints}. 
\end{proof}

\begin{remarks}\label{r4.1} 
There is a finer operation which consists of inverting $p$ on morphisms rather than on objects, but Corollary \ref{c5.2} is false for these categories. For example, the sheaf $\bigoplus_{n\ge 1} W_n$ is a non-zero object of $\RSC_\Nis$, but $h_0^{\rec}$ maps it to $0$ in $\HI_\Nis$. 
\end{remarks}

In the sequel, $(\ol{X},Y)$ is as in Theorem \ref{thm:h-h0}.

\subsection{Preliminary lemmas}\label{s5.1}
In the rest of this section,
we use a change of coordinates $\bcube \simeq (\P^1, 1)$
given by $\A^1 \to \P^1 \setminus \{ 1 \}, ~t \mapsto t/(t-1)$.
Let $\cube:=\bP^1-\{1\}$.
Take $S\in \Sm$ and a closed integral subscheme $V\subset S\times \cube \times X$
that is finite and surjective over $S \times \cube$.
We have a commutative diagram
\begin{equation}\label{eq1.proof-mainthm}
\xymatrix{ 
& V\ar@{^{(}->}[r]\ar[d] & X\times \cube\times S\ar[d]\\
\Vb^N \ar[r]\ar[d] & \Vb \ar@{^{(}->}[r]\ar[d] & \Xb\times \P^1\times S\ar[d]^{p}\\
\Wb^N \ar[r]\ar[rrd]_{\gamma} & \ol{W} \ar@{^{(}->}[r]\ar[r]& \Xb\times S\ar[d]\\
&& \Xb \\
}\end{equation}
where $\Vb$ is the closure of $V$ in $\Xb\times \bP^1\times S$, $\Wb$ is the image of $\Vb$ under the projection
$p$, and $\Vb^N\to \Vb$ and $\Wb^N\to \Wb$ are the normalizations. 
Let $\phi_V:\Vb^N\to \Xb\times \bP^1\times S$ be the natural map.
Let $\iota_\infty : \Xb\times S \to \Xb\times\bP^1\times S$ be induce by $\infty\in \bP^1$. Put 
\[
\partial^\infty\Vb=\iota_\infty^{-1}(\Vb)=p(\Vb\cap (\Xb\times\{\infty\}\times S)) \subset \Xb\times S.
\]
Putting $\Wb^o=\Wb\backslash \partial^\infty V$ and $\WbNo=\Wb^N\times_{\Wb} \Wb^o$, we have
\begin{equation}\label{eq:emb-into-line}
\Vb \times_{\Wb} \WbNo \subset \WbNo\times(\bP^1-\{\infty\}).
\end{equation}
Let $\ol{V}^o$ be the reduced part of an irreducible component 
of $\Vb \times_{\Wb} \WbNo$ which dominates $\WbNo$.
(Thus $\ol{V}^o \to \ol{V}$ is birational.)

\begin{lemma}\label{suslin-lemma-1}
If $\Wb^o=\emptyset$, then $V=W \times \cube$ with $W=\Wb\cap (X\times S)$.
\end{lemma}
\begin{proof}
The assumption implies $\Wb\subset p(\Vb\cap(\Xb\times\{\infty\} \times S))$ and hence
\[
\dim \Wb \leq \dim \Vb\cap(\Xb\times\{\infty\} \times S) <\dim \Vb.
\]
Noting $\Vb\hookrightarrow \Wb\times\bP^1$, we get $\Vb=\Wb\times \bP^1$, which implies the desired assertion.
\end{proof}

\begin{lemma}\label{suslin-lemma0}
If $\Wb^o\not=\emptyset$, $\Vb^o$ is finite over $\WbNo$.
\end{lemma}
\begin{proof}
$\Vb$ is proper over $\Wb$ so that $\Vb^o$ is proper over $\WbNo$.
On the other hand $\WbNo\times(\bP^1-\{\infty\})$ is affine over $\WbNo$ and so is $\Vb^o$. 
This implies the lemma.
\end{proof}

Now we consider the modulus condition for $V$:
\begin{equation}\label{modulus.eq}
\phi_V^{-1}(Y\times \bP^1\times S)\leq \phi_V^{-1}(\Xb\times \{1\} \times S).
\end{equation}
Let $y$ be the standard coordinate on $\bP^1-\{\infty\}=\Spec(k[y])$.
(Note that the divisor 
involved in the modulus condition
is $\{ 1 \} \subset \P^1 - \{ \infty \}$
defined by the ideal $(1-y) \subset k[y]$.)
Let
$I\subset \sO_{\WbNo}$ be the ideal sheaf of $Y\times_\Xb \WbNo \subset \WbNo$.

\begin{lemma}\label{suslin-keylemma}
Assuming $\Wb^o\not=\emptyset$, \eqref{modulus.eq} is equivalent to the conditions:
\begin{itemize}
\item[$(i)$]
$\Vb\cap (Y\times \cube \times S)=\emptyset$. 
\item[$(ii)$]
Locally on $\WbNo$, $\Vb^o$ is defined by an equation
\[
f(y):=(1-y)^m+ \underset{1\leq \nu \leq m}{\sum}\; a_\nu (1-y)^{m-\nu}\;\;\text{ with }\;
a_\nu\in \Gamma(\WbNo,I^\nu), 
\]
in $\WbNo\times(\P^1 - \{ 1 \})=\WbNo\times\Spec(k[y])$
(see \eqref{eq:emb-into-line}).
\end{itemize}
\end{lemma}

\begin{proof}
By Lemma \ref{suslin-lemma0}, the minimal polynomial over $k(\Wb)$
of the image of $y$ in $\Gamma(\Vb^o,\sO)$:
\[
f(t)= (1-t)^m+ \underset{1\leq \nu \leq m}{\sum}\; a_\nu (1-t)^{m-\nu} 
\]
has its coefficients $a_\nu\in A:=\Gamma(\WbNo,\sO)$.
We claim that $\Vb^o$ coincides with the closed subscheme $T\subset \WbNo\times\Spec(k[y])$ defined by the equation
$f(y) \in A[y]$. 
Indeed it is clear that $\Vb^o$ is contained in $T$,
hence it suffices to show that $T$ is integral.
Note that $T$ is a Cartier divisor in $\WbNo\times(\P^1 - \{ 1 \})$
which is finite over $\overline{W}^{N,o}$.
It follows that each irreducible component
dominates $\overline{W}^{N,o}$.
Hence the integrality is checked over the generic point, 
which holds by the irreducibility of $f$.
The claim is proved.
Thus we are reduced to showing the following. 

\begin{claim}\label{claim.suslin-keylemma}
The condition \eqref{modulus.eq} holds if and only if 
$\Vb\cap (Y\times \cube \times S)=\emptyset$ and $a_\nu\in \Gamma(\WbNo,I^\nu)$ for all $\nu$.
\end{claim}

The question is Zariski local and we may assume that $I$ is generated by $\pi \in \Gamma(\WbNo,\sO)$. 
Then \eqref{modulus.eq} holds if and only if $\ol{V}\cap (Y\times \cube \times S)=\emptyset$ and 
\begin{equation}\label{eq2-suslin-keylemma}
\theta:=\displaystyle{\frac{1-\overline{y}}{\pi}}\in \Gamma(\Vb^N\times_{\Wb^N} \WbNo,\sO).
\end{equation} 
Noting $\pi\in k(\Wb)$, the minimal polynomial of $\theta$ over $k(\Wb)$ is 
$$
g(t) =t^m+ \underset{1\leq \nu \leq m}{\sum}\; \frac{a_\nu}{\pi^\nu} \;t^{m-\nu}.
$$
Since $\Vb^o$ is finite over $\WbNo$ as is shown before, 
$\Vb^N\times_{\Wb^N} \WbNo$ is finite over $\WbNo$. Hence \eqref{eq2-suslin-keylemma} is equivalent to the condition that 
$\theta$ is integral over $\Gamma(\WbNo,\sO)$, which is equivalent to  
\[
\frac{a_\nu}{\pi^\nu}\in\Gamma(\WbNo,\sO)\quad\text{ for all }\nu.
\]
This proves the claim and the proof of Lemma \ref{suslin-keylemma} is completed.
\end{proof}

\subsection{Proof of Theorem \ref{thm:h-h0}}\label{subsect:screc-rec}
We put
\[ C_1(\ol{X}|Y):=\omega_!\uHom_{\MPST}(\Z_\tr(\bcube), \Z_\tr(M)) \in \PST
\]
and write by $\partial$ for the boundary map
$\delta_{1,0}^{0*} - \delta_{1,\infty}^{0*} :
C_1(\ol{X}|Y) \to \Z_\tr(M)$.
Fix $S\in \Sm$.
By Definitions \ref{d2.3-2} and \ref{def:old-rec}, 
it suffices to construct a homomorphism :
\begin{equation}\label{eq-1.proof-suslin}
\xi: \CXY 1 S \to \Phi(\Xb, Y)(S)
\end{equation}
such that the following diagram commutes:
\begin{equation}\label{eq0.proof-suslin}
\begin{CD}
 \CXY 1 S @>{\partial}>> \Ztr(X)(S) \\
 @VV{\xi}V @| \\
 \Phi(\Xb, Y)(S) @>{\tau}>>  \Ztr(X)(S)  \\ 
\end{CD}
\end{equation}
and such that we have
\begin{equation}\label{eq0-1.proof-suslin}
\Image(\tau)=\Image(\tau\circ\xi).
\end{equation}
Take a closed integral subscheme $V\subset S\times \cube \times X$, finite and surjective over $S\times \square$ and satisfying \eqref{modulus.eq}.
Consider the commutative diagram \eqref{eq1.proof-mainthm} and 
let $\phi: \Wb^N \to \Xb\times S$ be the induced map. 
We first suppose $\Wb^o\not=\emptyset$.
Then we have (see \S \ref{ssec:ref} for notations)
\[
(\Wb^N \rmapo{\phi} \Xb\times S)\;\in \sC_{(\Xb, Y)}(S).
\]
The projection $V\to \cube=\bP^1-\{1\}$ induces a rational function $g_V\in k(\Vb)^\times$.
By \cite[Prop.1.4 and \S 1.6]{Ful} we have
\begin{equation}\label{eq1.proof-suslin}
\partial V = \phi_* \div_{\Wb^N}(N g_V)\;\in \Ztr(X)(S) = \Cor(S, X),
\end{equation}
where $N:k(\Vb)^\times \to k(\Wb)^\times$ is the norm map induced by $\Vb\to \Wb$.
By Lemma \ref{suslin-keylemma} we have
\[
N g_V = f(0) = 1 + \underset{1\leq \nu \leq m}{\sum}\; a_\nu\; \in 
\Gamma(\WbNo,I)\subset G(\Wb^N,\gamma_\phi^*Y)
\subset \Phi(\Xb,Y)(S).
\]
We now define a map
\begin{equation}\label{eq1-1.proof-suslin}
\xi:\CXY 1 S \to \Phi(\Xb,Y)(S)
\end{equation}
by declaring
\[
\xi(V)= 
\begin{cases}
N g_V 
& \text{if} ~\Wb^o \not= \emptyset;
\\
0 & \text{if} ~\Wb^o=\emptyset.
\end{cases}
\]
Note that if $\Wb^o=\emptyset$,
then we have $\pd(V)=0$
by Lemma \ref{suslin-lemma-1}.
It follows that the diagram \eqref{eq0.proof-suslin} commutes
thanks to \eqref{eq1.proof-suslin}.

It remains to show \eqref{eq0-1.proof-suslin}.
To this end,
we take $(\phi_0 : \ol{C} \to \Xb \times S) \in \sC_{(\Xb,Y)}(S)$
and show
$\tau(G(\ol{C}, \gamma_{\phi_0}^*Y))
\subset \Image(\tau \circ \xi)$
(see \S \ref{ssec:ref} for notations).
Let $\Wb\hookrightarrow\Xb\times S$ be the image of $\phi_0$ and 
let $\Wb^N \to \Wb$ be its normalization
so that $(\phi : \Wb^N \to \Xb\times S)\in \sC_{(\Xb,Y)}(S)$.
Since 
$\tau(G(\ol{C}, \gamma_{\phi_0}^*Y))
\subset \tau(G(\Wb^N, \gamma_{\phi}^*Y))$,
it suffices to show the following.

\begin{lemma}\label{lem3-susulin}
%
%
The subgroup
$G(\Wb^N,\gamma_\phi^*Y)\subset \Phi(\Xb,Y)(S)$
is contained in the image of $\xi:\CXY 1 S \to \Phi(\Xb,Y)(S)$.
\end{lemma} 
\begin{proof}
Take  $g\in G(\Wb^N,\gamma_\phi^*Y)$.
Let $\Sigma\subset \Wb^N$ be the closure of the union of points $x\in \Wb^N$ of codimension one such that $v_x(g)<0$,
where $v_x$ is the valuation associated to $x$. 
Since $\Wb^N$ is normal, we have $g\in \Gamma(\Wb^N-\Sigma,\sO)$
and $g\in G(\Wb^N,\gamma_\phi^*Y)$ implies
\begin{equation}\label{eq4.proof-suslin}
g-1 \in \Gamma(\Wb^N-\Sigma,I),
\end{equation}
where $I\subset \sO_{\Wb^N}$ is the ideal sheaf of $\gamma_\phi^*Y\subset \Wb^N$.
Let 
\[
\psi_g: \Wb^N-\Sigma \to \bP^1-\{\infty\}
\]
be the morphism induced by $g$ and $\Gamma\subset \Wb^N \times \bP^1$
be the closure of the graph of $\psi_g$.
Let 
\[
\Vb\subset \Wb\times \bP^1\subset \Xb\times\bP^1 \times S
\]
be the image of $\Gamma$ under $\Wb^N \times \bP^1\to \Wb \times \bP^1$.
By \eqref{eq4.proof-suslin} we have $|\gamma_\phi^*Y|\subset \psi_g^{-1}(1)$ and hence
\begin{equation}\label{eq5.proof-suslin}
\Vb\cap (Y\times \cube \times S)=\emptyset
\end{equation}
so that 
\[
V:=\Vb\cap (\Xb\times \cube \times S)\subset X\times \cube\times S.
\] 
It suffices to show the following.

\enlargethispage*{20pt}

\begin{claim}\label{claim3.susulin}
$V\in \CXY 1 S$ and $\xi(V)=g$.
\end{claim}

Once we prove the first assertion, the second follows easily from the construction of $\xi$.
To prove the first assertion,  by \eqref{eq5.proof-suslin}, the map $V\to \cube\times S$ is proper and hence finite since $X$ is quasi-affine by the assumption. 
Moreover it is surjective since $\dim V =\dim \Wb =\dim S +1$. 
Hence it suffices to check the condition $(ii)$ of Lemma \ref{suslin-keylemma}. By definition 
\begin{enumerate}
\item[$(\spadesuit)$]
$\Gamma\cap \big((\Wb^N-\Sigma)\times (\bP^1-\{\infty\})\big)$
is the graph of $\psi_g$ and hence is defined by $y-g$ where $y$ is the standard coordinate of 
$\bP^1-\{\infty\}=\Spec(k[y])$.
\end{enumerate}
We have a diagram of schemes
\begin{equation}\label{eq5.suslin}
\xymatrix{
\Wb \ar[r]^{\hskip -20pt\iota_\infty} & \Wb\times\bP^1 & \ar[l] \Vb 
\\
\Wb^N \ar[r]^{\hskip -20pt\iota_\infty}\ar[u]^{\pi_{\Wb}} \ar[rd]_{id_{\Wb}} 
& \Wb^N\times\bP^1 \ar[u]\ar[d]^{pr} 
& \ar[l]\ar[u]\ar[ld] \Gamma' :=\Vb\times_{\Wb} \Wb^N
& \ar[l]\ar[lu]_{\pi_{\Vb}}\ar[lld]^{pr_\Gamma} \Gamma 
\\
&\Wb^N
}
\end{equation}
where $\iota_\infty$ are induced by $\infty\in \P^1$.
The natural map $\Gamma\to \Gamma'$ 
is a closed immersion onto an irreducible component
that dominates $\ol{V}$.
We claim
\begin{equation}\label{eq6.susulin}
\Sigma \subset 
\iota_\infty^{-1}(\Gamma).
\end{equation}
The claim implies 
$\WbNo:=\Wb^N\times_{\Wb}(\Wb\backslash \iota_\infty^{-1}(\Vb)) \subset\Wb^N-\Sigma$. 
Let
$\Vb^o$ be the reduced part of the irreducible component of
\[
{\Vb^o}':=\Vb\times_{\Wb}\WbNo =\Gamma'\times_{\Wb^N}\WbNo \subset \WbNo\times(\bP^1-\{\infty\})
\]
which dominates $\ol{V}$;
see the following diagram:
\[
\xymatrix{
&
\WbNo \ar@{^{(}->}[d] \ar@{^{(}->}[dl]
&
\WbNo \times(\P^1 - \{ \infty \}) \ar@{^{(}->}[d] \ar[l]
& 
{\Vb^o}' \ar@{^{(}->}[d] \ar@{_{(}->}[l]
& \Vb^o \ar@{^{(}->}[d] \ar@{_{(}->}[l]
\\
\Wb^N - \Sigma \ar@{^{(}->}[r]
&
\Wb^N \ar[d]
&
\Wb^N \times \P^1 \ar[d] \ar[l]
& \Gamma' \ar[d] \ar@{_{(}->}[l]
& \Gamma \ar[dl] \ar@{_{(}->}[l]
\\
&
\Wb
&
\Wb \times \P^1 \ar[l]
& \ol{V}. \ar@{_{(}->}[l]
}
\]
By $(\spadesuit)$, $\Vb^o$ is defined in $\WbNo\times \Spec(k[y])$ by the equation $y-g$ and thus
$V$ satisfies Lemma \ref{suslin-keylemma} $(ii)$.

It remains to show \eqref{eq6.susulin}. From \eqref{eq5.suslin}, it is equivalent to 
\begin{equation}\label{eq7.susulin}
\Sigma \subset pr\big((\Wb^N\times {\infty})\cap \Gamma\big).
\end{equation}
Since $pr_\Gamma$ is proper birational and $\Wb^N$ is normal,
$pr_\Gamma$ is an isomorphism above all codimension one points in $\Wb^N$ (but not necessarily in all codimension one points of $\Gamma$).
For a generic point $x\in \Sigma$, there is a unique codimension one point $y\in \Gamma$ such that $x=pr_{\Gamma}(y)$ 
and we have $v_y(g)=v_x(g)<0$ for $g\in k(\Gamma)=k(\Wb^N)$.
The projection $\Wb^N\times\bP^1\to \bP^1$ induces a morphism
$\Gamma \setminus (\Wb^N\times \{\infty\}) \to \bP^1-\{\infty\}$,
which corresponds to $g$. Hence we must have $y\in (\Wb^N\times \{\infty\})\cap \Gamma$
which proves \eqref{eq7.susulin}
by the properness of $pr_\Gamma$.
This completes the proof of Lemma \ref{lem3-susulin}.
\end{proof}

\enlargethispage*{30pt}

\enlargethispage*{20pt}

\end{document}